\documentclass[12pt]{amsart}

\usepackage{amssymb,graphicx,mathtools,amsthm, float}
\newtheorem{theorem}{Theorem}[section]

\newcommand{\pZ}{\mathbb{Z}_{\geq 0}}
\newcommand{\hF}{\hat{F}}
\title{Dehn's Algorithm for Simple Diagrams}

\author{Charles Frohman }
\address{ Department of Mathematics, The University of Iowa}
\email{charles-frohman@uiowa.edu}
\author{Joanna Kania-Bartoszynska}
\address{ Division of Mathematical Sciences, The National Science Foundation}
\email{jkaniab@nsf.gov}
\thanks{This material is based upon work supported by and while serving at the National Science Foundation. Any opinion, findings, and conclusions or recommendations expressed in this material are those of the authors and do not necessarily reflect the views of the National Science Foundation.}

\begin{document}
\begin{abstract} A version of Dehn's algorithm for simple diagrams on a once punctured surface representing simple diagrams on a closed surface is presented.\end{abstract}

\maketitle

 \section{Introduction}
 
Let $F$ be a closed surface of genus $g>0$.  The  fundamental group of $F$ can be presented as
\[ \pi_1(F)=\{a_i,b_i| \prod_{i=1}^g a_ib_ia_i^{-1}b_i^{-1}=e\}.\]
Max Dehn \cite{D2} gave  algorithms for deciding if a word in the $a_i^{\pm 1}$ and $b_i^{\pm 1}$ represents the identity, and to decide if two words represent the same conjugacy class.  In this note we extend his algorithm for conjugacy to simple diagrams on a once punctured surface of negative Euler characteristic representing diagrams on the closed surface obtained by filling in the puncture.

Dehn's algorithm for conjugacy classes proceeds as follows.  If $w$ is a cyclic word, begin by eliminating all occurrences of a letter adjacent to its inverse. The result is {\bf freely cyclically reduced}.  If the greater half of the relator or its inverse, treated as cyclic words, appears in the cyclic word, replace it by the inverse of the shorter half.  Continue these steps until there are no long halves of relators appearing in the cyclic word  and the word is freely cyclically reduced. The resulting word is a least length representative of the given conjugacy class. If the word contains no half relators, then the cyclic word is the unique least length representative of the conjugacy class.  If it has half relators, then there are finitely many least length representatives coming from swapping half relators. Two cyclic words are conjugate if and only if they have the same collection of least length representatives.

Suppose that $\hat{F}$ is a closed surface of nonpositive Euler characteristic  and $p\in \hat{F}$. Let $F=\hat{F}-\{p\}$. Choose an ideal triangulation of $F$ with edges $C$.  The ends of all the edges of the triangulation
are mapped to the end of $F$. This gives rise to a cyclic ordering of the ends of the edges in $C$  and a cyclic ordering of all the corners of the triangles. We say an edge of the triangulation is {\bf antipodal} if in the cyclic ordering there are the same number of corners of triangles between the two ends of the edge on either side in the cyclic ordering. We define {\bf antipodal corners} similarly.

A {\bf simple diagram} on a surface is a disjoint collection of simple closed curves, none of which bounds a disk on the surface.
The inclusion map $F\rightarrow \hat{F}$ sends a simple diagram on $F$ to a disjoint union of simple closed curves on $\hat{F}$. The system of curves on $\hat{F}$ is not necessarily a simple diagram because any simple closed curve on $F$ that bounds a punctured disk on $F$  bounds a disk on $\hF$.
Conversely any simple diagram on $\hat{F}$ can be perturbed so that it misses the puncture, giving rise to  a simple diagram on $F$.  The simple diagram on $F$ is not unique as different perturbations differ by sliding the diagram over the puncture.
If the image under inclusion  of a simple diagram $S'$ on $F$ into $\hat{F}$, is isotopic to the simple diagram $S$ in $\hat{F}$, then we say that $S'$ {\bf represents} $S$.

If the set of lines $C$ defines an ideal triangulation of $F$, and $S'$ is a simple diagram on $F$, then there is a  function $f_{S'}:C\rightarrow \pZ$ that assigns the geometric intersection number $i(c,S')$ of $S'$ to each edge $c\in C$.  This function has the properties that, for any three edges $a$, $b$, $c$ that bound a triangle, the sum $f_{S'}(a)+f_{S'}(b)+f_{S'}(c)$ is even, and the integers $f_{S'}(a),f_{S'}(b)$ and $f_{S'}(c)$ satisfy the triangle inequality. The functions $f:C\rightarrow \pZ$ satisfying these conditions are called {\bf admissible colorings } of $C$, and are in one to one correspondence with isotopy classes of simple diagrams on $F$.

In this paper we give an algorithm to decide if two simple diagrams on $F$ represent the same diagram on $\hat{F}$.

The {\bf weight} of a simple diagram corresponding to the admissible coloring $f_{S'}:C\rightarrow \pZ$ is $\sum_{c\in C}f_{S'}(c)$.  

\vspace{.1in}

{\bf Theorem \ref{seq}} {\it Suppose that $S'$ is a simple diagram on the punctured surface $F$ that represents a simple diagram $S$ on the closed surface $\hF$. Suppose that an ideal triangulation of $F$ has been chosen.  There is a sequence of simple diagrams $S'=S'_0,S'_1,\ldots, S'_n$ on $F$ all representing $S$ such that $S'_{i+1}$ is obtained by a handleslide across the puncture from $S'_i$,  the weights of the diagrams are monotonically decreasing, and the final diagram $S'_n$ is a least weight representative of $S$.  }

\vspace{.1in}

If $\{a,b,c\}$ are the edges of an ideal triangle, the triangle has  three ``corners'' where two edges of the triangle meet.  If $f:C\rightarrow \pZ$ is an admissible coloring of an ideal triangulation  then the corner number associated to the corner between $a$ and $b$ is \[ \frac{f(a)+f(b)-f(c)}{2},\] with the remaining corner numbers defined analogously.  The corner numbers tell you how many strands of the diagram run between the two sides at the corner, when the diagram realizes its geometric intersection numbers with the edges of the triangulation.

The corner numbers of the triangulation give rise to {\bf bands}, where in the cyclic ordering on the corners, consecutive strings of corner numbers are positive, and {\bf gaps} where the corner numbers are zero.  A band is a {\bf half band} if the gaps at either end are antipodal in the cyclic ordering.  A band is {\bf long} if it runs more than halfway around the vertex in the cyclic ordering.

\vspace{.1in}

{\bf Theorem \ref{lite}}  {\it There is a unique weight minimizing diagram representing the simple diagram $S\in \hF$ unless the weight minimizer has a maximal half band, or  there is an antipodal edge in the triangulation and the diagram contains a component isotopic to a simple closed curve corresponding to that edge. }  

\vspace{.1in}

A variant of the second clause involving antipodal edges also shows up in the extension of Dehn's algorithm to one relator presentations of a surface group that have a letter appearing  antipodally in the cyclic word of the relator.
If there is more than one weight minimizing diagram then all the weight minimizing diagrams are related by sliding half bands across the puncture,
or sliding push offs of an antipodal edge across an annulus containing the puncture.

\section{Curves and lines on finite type surfaces}

We set up the combinatorial framework for describing and proving Dehn's algorithm for simple diagrams on a closed oriented surface.

A surface $F$ has finite type if there is a closed oriented surface $\hF$ and a finite collection of points $\{p_i\}\subset \hF$ so that $\hF-\{p_i\}=F$.   The points $\{p_i\}$ are called punctures.

Suppose that $X$ and $Z$ are properly embedded $1$-manifolds in a finite type surface $F$ where  $X$ is compact.   We say that $X'$ is a transverse representative of $X$ with respect to $Z$, if $X'$ is ambiently isotopic to $X$ via a compactly supported isotopy, and $X'\pitchfork Z$.   Define the {\bf geometric intersection number} of $X$ and $Z$, denoted $i(X,Z)$, to be the minimum cardinality of $X'\cap Z$ over all transverse representatives of $X$. We could have instead worked with $Z$ up to compactly supported ambient isotopy and taken the minimum over all $Z'$ isotopic to $Z$ and transverse to $X$ and obtained the same number, thus $i(X,Z)=i(Z,X)$.  

 A {\bf  bigon} with respect to $1$-manifolds $X$ and $Z$ embedded in $F$ is a disk $D\subset F$ such that the boundary of $D$ consists of the union of two arcs $a\subset X$ and $b\subset Z$. If there is a bigon, there is always an innermost bigon, whose interior is disjoint from $X\cup Z$.
A transverse representative of $X$ realizes the geometric intersection number $i(X,Z)$ if and only if there are no  bigons,  \cite{FLP}.

 It is always possible to put a complete Riemannian metric on $F$, so we assume we have one.
If $C$ is a properly embedded disjoint system of lines in $F$ define a metric on the components $D$ of $F-C$ by
\[ d(p,q)=\inf{\{length(\alpha)|\alpha:[0,1]\rightarrow D \ \mathrm{is \ smooth}, \alpha(0)=p,\alpha(1)=q}\}.\] 
If $D$ is a component of $F-C$ let $\Delta$ be the metric space completion of $D$.
The surface $\Delta$ is diffeomorphic to the result of removing finitely many points from a compact surface.  The boundary of $\Delta$ is a finite collection of lines.  If $\Delta$ is diffeomorphic to a disk with three points removed from its boundary then $\Delta$ is an {\bf ideal triangle}.    

If $C$ is a disjoint system of properly embedded lines in $F$ and the completion of every component of the complement of $C$ is an ideal triangle we say that $C$ defines an {\bf ideal triangulation} of $F$. Alternately, we could define an ideal triangulation as the result of identifying a collection of ideal triangles along boundary components in pairs, along with a homeomorphism to $F$.  A noncompact finite type surface admits an ideal triangulation if and only if its Euler characteristic $\chi(F)$  is negative.  An ideal triangulation of $F$ is made up of $-2\chi(F)$ ideal triangles sharing $-3\chi(F)$ edges. 

\begin{figure}[htbp] \begin{center}\scalebox{1}{\includegraphics{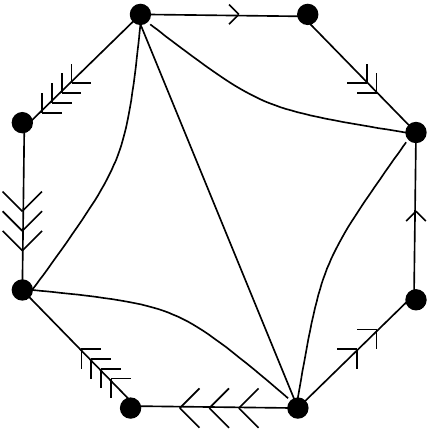}}\end{center}\caption{An ideal triangulation of $\Sigma_{2,1}$}\label{gen}\end{figure}

In Figure \ref{gen} we depict  an ideal triangulation of $\Sigma_{2,1}$ the once punctured surface of genus $2$.
The sides of the octagon are identified as indicated. The sides along with the lines in the interior of the octagon define an ideal triangulation with $6$ triangles and $9$ edges.

There are two kinds of ideal triangles in an ideal triangulation, {\bf embedded ideal triangles} and {\bf folded ideal triangles}. Embedded ideal triangles are just that, folded ideal triangles are triangles with two edges identified.  Hence a folded ideal triangle only has two edges, one of which has multiplicity $2$.  If the surface $F$ only has one puncture there are no folded triangles, as a folded triangle abuts two distinct punctures.

You can pass from any ideal triangulation to any other by a sequence of $1-1$ Pachner moves.  The $1-1$ Pachner move applies to two ideal triangles that share an edge. The two triangles form a quadrilateral, where the shared edge is a diagonal. The Pachner move replaces the shared edge by the other diagonal.

\begin{figure}[htbp]\begin{center} \includegraphics{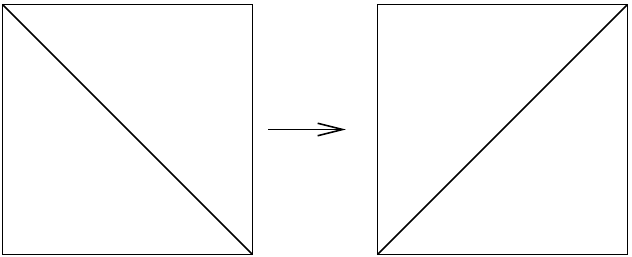} \end{center}\caption{The $1-1$-Pachner move.}\end{figure}

Suppose that $C$ is the set of edges of an ideal triangulation with no folded triangles,  and  $f:C\rightarrow \pZ$ is an admissible coloring.
If $\{a,b,c\}$ are the sides of an embedded triangle then  there is a collection of arcs properly embedded in the triangle so that no arc has both endpoints in the same edge and there are $\frac{f(a)+f(b)-f(c)}{2}$ arcs running between the sides $a$ and $b$, $\frac{f(a)+f(c)-f(b)}{2}$ arcs running between the sides $a$ and $c$, and $\frac{f(b)+f(c)-f(a)}{2}$ arcs running between the sides $b$ and $c$.  The pairs of sides are the {\bf corners} of the triangle and the numbers of strands are the {\bf corner numbers } of $f$ associated to the given corners.

On the other hand if $S$ is a simple diagram, it can be isotoped so that $S\pitchfork C$ and there are no bigons  between $S$ and $C$.  Inside the triangles the components of $S$ are arcs running between the sides, with no bigons.  The function $f:C\rightarrow \pZ$ given by $f(c)=i(S,c)$ is an admissible coloring of $C$.

Putting the pieces together;

\begin{theorem} Given a surface $F$ of finite type with an ideal triangulation determined by the set of properly embedded lines $C$,  there is one-to-one correspondence between simple diagrams on  $F$ up to isotopy and admissible colorings of $C$. \end{theorem}
\qed

Now suppose that $\hF$ is closed, $p\in \hF$ is a point.  The surface $F=\hF-\{p\}$ has a single end. Therefore every end of every edge, and every corner of every triangle in an ideal triangulation abuts the end where the point was deleted.   This gives rise to cyclic orderings of the ends of the edges and of the corners of the triangles.

\begin{figure}[htbp]\begin{center}\includegraphics{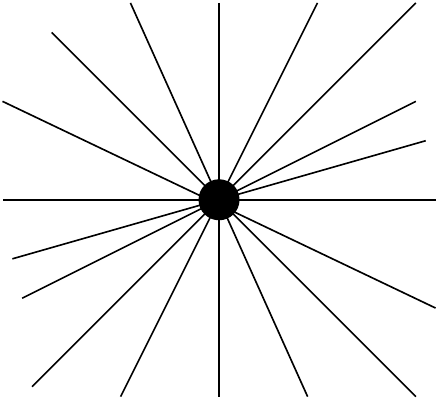}\caption{The ends of the edges at the puncture.}\end{center}\end{figure}

A {\bf band} $A_1,A_2,\ldots,A_k$ of an admissible coloring $f:C\rightarrow \pZ$
of an ideal triangulation of a once punctured surface $F$ is an ordered family of corners of triangles so that they are consecutive going counterclockwise in the cyclic ordering of corners at the end of the surface, and all of the corner numbers of $f$ at the $A_i$  are positive.  

Suppose that $A_0$ is the corner immediately preceding $A_1$ and $A_{k+1}$ is the corner immediately following $A_k$. We say the band is {\bf maximal} if the corner numbers of $A_0$, and $A_{k+1}$ with respect to $f$ are zero.  If the corner number is zero we say that corner is a {\bf gap}.  If the band has more than half of all the corners, we say it is a {\bf long band}. If the gaps at either end of a  maximal band are antipodal  we call it a {\bf half band}. If a maximal band has less than half the corners we call it a {\bf short band}. 

\begin{figure}[htpb]\begin{center}\includegraphics{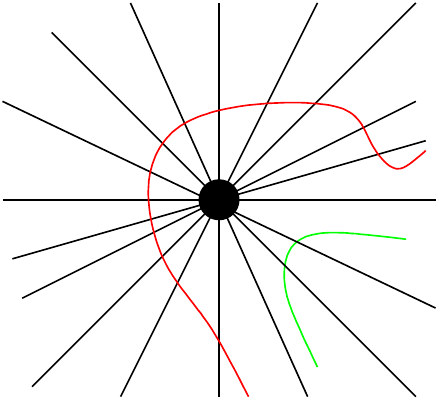} \end{center}\caption{Two bands.}\label{redband}\end{figure}

Given a band, there is a strand of the corresponding diagram running close to the vertex through the corners of the band. Close means that no other part of the diagram touches the ends of the edges heading towards the corners of the band.  The red strand in Figure \ref{redband} corresponds to a long band, and the green strand to a short band. Each triangle has three sides, of which two are involved in any corner. We are assuming the bands  leave the triangles that contain their ends in the  side of the triangle not belonging to that corner. There are four gaps in Figure \ref{redband}.

\section{Dehn's algorithm for simple diagrams}
Suppose as before  that $\hF$ is closed, $F=\hF-\{p\}$, and $S$ is a simple diagram on $\hF$.  After a small perturbation, $S\cap \{p\}=\emptyset$.
Hence there are simple diagrams $S'$ on $F$ that under the inclusion $F\hookrightarrow \hF$ are isotopic to $S$. Recall that in this case we say that $S'$ represents $S$.  Recall also that the  diagram $S'$ is a  weight minimizing  representative of $S$ if in addition its weight $\sum_{c\in C}f_{S'}(c)$  is minimal among all representatives.

A {\bf handleslide} between two representatives $S'$ and $S''$ is an isotopy in $\hF$ that moves generically across the point $p$. Our goal is to describe an algorithm that takes a representative of $S$ to a minimal weight representative though a sequence of handleslides, such that the  weights of the intermediate diagrams are monotonically decreasing. 

There are two kinds of handleslides that are particularly simple. One of them involves an annulus  and the other a bigon, in each case with $p$ in its interior. 
We recall the {\bf Annulus theorem}:
\begin{theorem}\cite{E}\label{Annulus}
Suppose that $J$ and $J'$ are two disjoint simple closed curves on a surface $F$ that do not bound disks. If $J$ is isotopic to $J'$ then there exists an annulus embedded in $F$ whose boundary is $J\cup J'$.  
\end{theorem}
If $S$ and $S'$ are two disjoint non-isotopic diagrams on $F$ representing the same simple diagram on $\hF$ then by \ref{Annulus} there are annuli in $\hF$ cobounding the components of $S$ and $S'$. Since $S$ and $S'$ are not isotopic in $F$ one of these annuli must be punctured.  
Bigons arise when $S$ and $S'$  represent the same diagram in $\hF$ and
$i(S,S')\neq 0$. In this case there is a bigon on $\hF$ between  $S$ and $S'$ that gives rise to a punctured bigon on $F$.   We use these as guides for isotopies.

\begin{theorem} \label{seq} Suppose that $S'$ is a simple diagram on the punctured surface $F$ that represents a simple diagram $S$ on the closed surface $\hF$.  Suppose that an ideal triangulation $C$ of $F$ has been chosen.  There is a sequence of simple diagrams $S'=S'_0,S'_1,\ldots, S'_n$ on $F$,   all representing $S'$, such that $S_{i+1}'$ is obtained by a handleslide across the puncture from $S_i'$, the weights of the diagrams are monotonically decreasing, and the final diagram $S_n'$ is a least weight representative of $S'$.  \end{theorem}

\proof  Suppose that $S'$ and $S''$ are representatives of $S$, and that $S''$ has least weight among all simple diagrams on $F$ that represent $S$. We find a sequence of diagrams related by handleslides taking $S'$ to $S''$ whose weights are monotonically decreasing.    We break the argument into two cases.  The first case we consider is when $i(S',S'')=0$.  We can assume that if any components of $S'$ and $S''$ are isotopic on $F$ then we have pulled them close to one another, so that they cobound an annulus on $F$ and that annulus contains no other components of $S'$ or $S''$.  The complexity in this case is the number of components of $S'$ that have not been paired this way with components of $S''$.  Our moves always increase the number of paired components.  Once we have arrived at two diagrams that are completely paired, the final isotopy is obvious.  In the second case when $i(S',S'')\neq 0$ the complexity will be the geometric intersection number of the two diagrams. In this case, we can find an innermost bigon on $\hF$ that contains $p$. We will then use this bigon as a guide that will decrease the geometric intersection number of the diagrams being considered without increasing weight. After finitely many steps we will be in the first case.

We start with the case where $i(S',S'')=0$.  Assume we have built the sequence $S'=S_0'$, $S_1,S_2',\ldots, S_i'$ with monotonically decreasing weights.
Suppose that we have chosen representatives of the isotopy classes of $S_i'$ and $S''$  in $F$ so that $S_i'\pitchfork S''$, and $S_i'$ and $S''$ realize their geometric intersection numbers with each other and with the edges of the triangulation $C$. Assume we have a collection of paired components of $S_i'$ and $S''$, in the sense that they cobound annuli that contain no other components.  Since the images of $S_i'$ and $S''$ are isotopic in $\hF$,  but not in $F$ there must be an annulus in $\hF$  that contains $p$, whose boundary consists of a component of $S_i'$ and  a component of $S''$ that have not been paired.  This annulus might have nonempty intersection with one of the pairing annuli. In this case slide the pairing annulus off the punctured annulus. 

\begin{figure}[htbp]\begin{center}\includegraphics{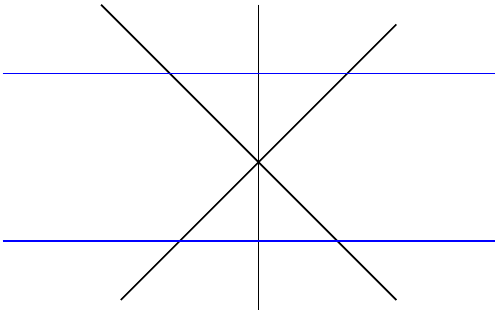}\end{center}\caption{A punctured annulus}\label{annul}\end{figure}

In Figure \ref{annul} we depict a punctured annulus.  In the diagram the two horizontal blue curves are the boundary components of the annulus.  In general one boundary component of the annulus will have greater geometric intersection number with $C$ than the other. The more weighty of the two components must be a component of $S'_i$. The annulus gives a guide for a handleslide to a representative $S'_{i+1}$ of lower weight, with more components of $S'_{i+1}$ paired with components of $S''$ than before.

In  Figure \ref{annul} the boundary components of the punctured annulus have the same weight with respect to $C$.
The regions to the left and right of the puncture are the same, so the two ends of edges coming diagonally on the right, and the two ends of edges coming diagonally from the left, are ends of edges of the same triangle. Since a triangle only has three sides, two of the ends of edges belong to the same edge of the triangle.  Since $S_i'$ and $S''$ realize their geometric intersection numbers with the edges of triangulation, it must be that the pair of ends of the same edge is either the upper one on the left and the lower one on the right or the upper one on the right and the lower one on the left.   This can only happen if there is an antipodal edge, and the components of $S_i'$ and $S''$ that cobound this annulus are isotopic to simple closed curves that are pushoffs of this edge. The annulus can still be used as a guide for a handleslide from $S_i'$ to $S_{i+1}'$ that has the same weight, and we have increased the number of paired components.

\begin{figure}[htbp]\begin{center}\includegraphics{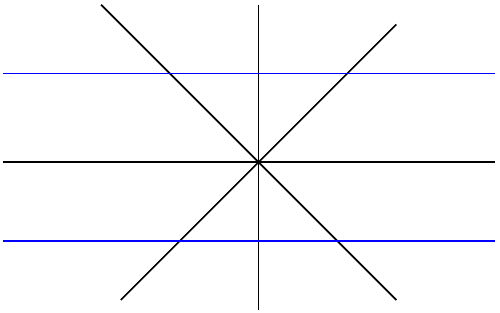}\end{center}\caption{An antipodal edge}\label{antipodal} \end{figure}

The punctured annulus could also look like the diagram in Figure \ref{antipodal}, where there is an edge of the triangulation lying completely inside the annulus. Note that Euler characteristic considerations keep more than one edge of the triangulation from being contained completely inside the annulus. The analysis is the same as for the last punctured annulus, only it is more obvious that for the boundary components of the annulus to have the same weight the triangulation must have an antipodal edge.

We repeat the last step until we have performed a sequence of handleslides that at the very least do not increase weight, arriving at the least weight $S'_n$ that has its components completely paired with the components of $S''$.

Assume we have built the sequence $S'=S_0'$, $S_1',S_2',\ldots, S_i'$ with monotonically decreasing weight, and decreasing geometric intersection number with $S''$.
Now suppose that $S_i'$ and $S''$ represent $S$, but $i(S'_i,S'')\neq 0$. Assume that $S_i'$ and $S''$ realize their geometric intersection number with $C$, are transverse to one another and realize their geometric intersection number. 
Finally assume that $S_i'\cap S''\cap C=\emptyset$.

Using the inclusion map $F\rightarrow \hF$, and the fact that the images of $S'_i$ and $S''$ are isotopic in $\hF$, there exists a bigon $B$ in $\hF$ whose boundary is  made up of an arc $a$ of $S'_i$ and an arc $b$ of $S''$ and whose  interior is disjoint from $S_i'\cup S''$.  It must be that $p$ is in the interior of $B$, otherwise $S'$ and $S''$ would not realize their geometric intersection number in $F$.

\begin{figure}[htbp]\begin{center}\includegraphics{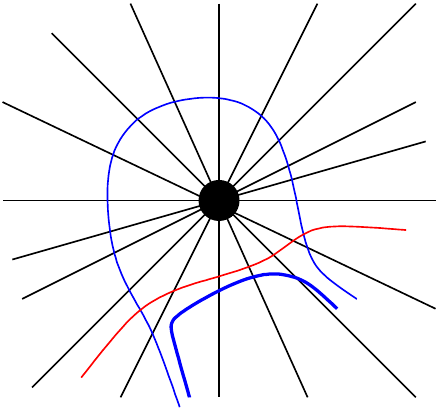} \end{center} \caption{Handleslide over a bigon} \label{longerside}\end{figure}

The sides of the bigon are made up of parts of bands.  If one of the parts of the band is longer than the other, then the longer side belongs to $S_i'$, see Figure \ref{longerside}.

\begin{figure}[htbp]\begin{center}\includegraphics{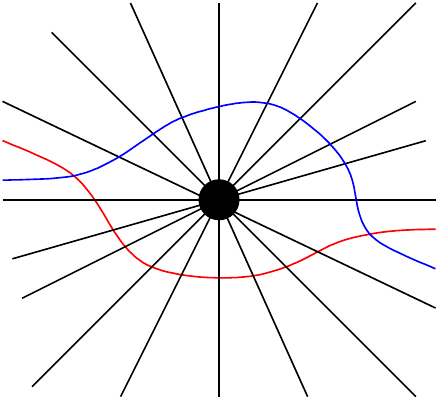}\end{center}\caption{ A half bigon}\label{halfbigon}  \end{figure}

The  corresponding handleslide allows us to pass from $S'_i$ to $S_{i+1}'$ which has lower weight, and has strictly lower geometric intersection number with  $S''$.  In Figure \ref{longerside} the handleslide replaces the longer blue strand with the shorter blue strand.
There is only one case where there is no shortening and that is when the bigon is made up of two maximal half bands like in Figure \ref{halfbigon}.

Slide  the arc of the  bigon belonging to $S'_i$ to obtain $S'_{i+1}$ without  increasing the weight and so that   $i(S_{i+1}',S'')<i(S_i',S'')$.

\qed

How many weight minimizing representatives are there?

\begin{theorem} \label{lite} There is a unique weight minimizing diagram representing the simple diagram $S\in \hF$ unless the weight minimizer has a maximal half band, or  there is an antipodal edge in the triangulation and the diagram contains a component isotopic to the simple closed curve corresponding to that edge.\end{theorem}

\proof  Suppose that there are two least weight diagrams $S'$ and $S''$ on $F$ representing the same simple diagram $S$ on $\hF$. By Theorem  \ref{seq} there are sequences of diagrams representing $S$ that are monotone decreasing in weight starting at  $S'$ and  terminating at $S''$.  Since $S'$ and $S''$ are both least weight this implies that in the  sequence of handleslides over bigons and punctured tori the weight of the diagrams is constant. As noted in  the proof of Theorem  \ref{seq} there are only two ways this can happen: via sliding by half bands, or if there is an antipodal edge.

 \qed

\end{document}